\documentclass[11pt]{amsart}

\usepackage{amsmath}
\usepackage{amstext}
\usepackage{amssymb}
\usepackage{amsthm}
\usepackage{enumerate}
\usepackage{mathrsfs}
\usepackage[all]{xy}
\usepackage{datetime}

\makeatletter
\def\imod#1{\allowbreak\mkern10mu({\operator@font mod}\,\,#1)}
\makeatother

\newtheorem{theorem}{Theorem}[section]

\newtheorem{prop}[theorem]{Proposition}
\newtheorem{lemma}[theorem]{Lemma}
\newtheorem{claim}[theorem]{Claim}
\newtheorem{corollary}[theorem]{Corollary}

\theoremstyle{definition}

\theoremstyle{remark}
\newtheorem{remark}[theorem]{Remark}

\theoremstyle{remark}

\numberwithin{equation}{section}

    \DeclareMathOperator{\diam}{diam}

    \DeclareMathOperator{\osc}{osc}
    
    \DeclareMathOperator{\IS}{IS}
    \DeclareMathOperator{\tr}{tr}

    \newcommand{\restrict}{\!\upharpoonright\!}

    \newcommand{\inx}[1]{\ \epsilon_{#1}\ }

\def\R{{\mathbb R}}

\def\N{{\mathbb N}}

\def\Q{{\mathbb Q}}

\title[$\mathbf{\Sigma}^1_2$ counterparts to statements that are equivalent to the Continuum Hypothesis]{the $ \mathbf{\Sigma}^1_2$ counterparts to statements that are equivalent to the Continuum Hypothesis}

\author{Asger T\"ornquist}

\address{Department of Mathematics, University of Copenhagen, Universitetsparken 5, 2100 Copenhagen, Denmark}

\email{asgert@math.ku.dk}

\author{William Weiss}

\address{Department of Mathematics, University of Toronto, 40 St. George St., Toronto, Ontario, Canada.}

\email{weiss@math.toronto.edu}

\thanks{Asger T\"ornquist was supported by a Sapere Aude fellowship (level 2) from DenmarkÕs Natural Sciences Research Council, no. 10-082689/FNU, and a Marie Curie re-integration grant, no. IRG-249167, from the European Union. William Weiss was supported in part by an NSERC discovery grant.}



\date{\today\  (\xxivtime)}

\begin{document}

\maketitle

\begin{abstract}
We consider natural $\Sigma^1_2$ definable analogues of many of the classical statements that have been shown to be equivalent to CH. It is shown that these $\Sigma^1_2$ analogues are equivalent to that all reals are constructible. We also prove two partition relations for $\Sigma^1_2$ colourings which hold precisely when there is a non-constructible real. 
\end{abstract}

\section{Introduction}

In the mathematical literature, one finds a great number of statements that have been proved to be equivalent to the Continuum Hypothesis (CH). One well-known such equivalence is due to Sierpinski, and states that CH is equivalent to that the plane $\R^2$ is the union of two sets $A,B\subseteq\R^2$ such that each horizontal section of $A$ is countable, and each vertical section of $B$ is countable. Another example is Davies' theorem, which states that CH is equivalent to that every function $f:\R^2\to\R$ admits a representation
$$
f(x,y)=\sum_{n=0}^\infty g_n(x)h_n(y),
$$
where $g_n,h_n:\R\to\R$ are functions and the sum above has only finitely many non-zero terms for every $(x,y)\in\R^2$.

In these types of theorems, usually the direct implication from CH is proved by a straight-forward inductive construction by well-ordering the reals in order type $\omega_1$, and exploiting that each initial segment is countable. The result of the construction will usually be definable from the well-ordering. Perhaps it is no surprise then that if we work in G\"odel's constructible universe $L$ where there is a canonical choice of a well-ordering of $\R$, which moreover is $\Sigma^1_2$, then with some care it can be shown in many cases that there are $\Sigma^1_2$ definable witnesses to the direct implication.

On the other hand, the reverse implication often requires considerable ingenuity and does not at first seem to conform to a set pattern. In light of the above discussion about the situation in $L$, it is natural to ask what happens if we take a statement which implies CH, and replace it with a corresponding $\Sigma^1_2$ version. In \cite{towe09} we considered the $\Sigma^1_2$ counterpart of Davies' theorem, and showed the following ``$\Sigma^1_2$ Davies' Theorem'': All reals are constructible $(\R\subseteq L$) if and only if every $\Sigma^1_2$ function $f:\R^2\to\R$ admits a representation
$$
f(x,y)=\sum_{n=0}^\infty g(x,n)h(x,n),
$$
where $g,h:\R\times\omega\to\R$ are $\Sigma^1_2$ functions, and the sum above has finitely many non-zero terms at each $(x,y)\in\R^2$. 

It is natural to ask if this type of definable converse, which was found in the case of Davies' theorem, could hold for some of the many other statements that are equivalent to CH. However, the proof in \cite{towe09} did not give a clear indication in this direction. In this paper we will prove that a number of the classical CH equivalents admit natural $\Sigma^1_2$ counterparts which turn out to be equivalent to that all reals are constructible. Specifically:

\begin{theorem}\label{t.mainthm1}
The following statements are equivalent:
\begin{enumerate}
\item $\R\subseteq L$.
\item There are $\Sigma^1_2$ sets $A,B\subseteq\R^2$ such that $A\cup B=\R^2$, and all the sections $A_x=\{y\in \R: (x,y)\in A\}$ and $B^y=\{x\in \R: (x,y)\in B\}$ are countable.
\item There are $\Sigma^1_2$ sets $A_1,A_2,A_3\subseteq\R^3$ such that $A_1\cup A_2\cup A_3=\R^3$, and every line $l$ in the direction of the $x_i$-axis meets $A_i$ in finitely many points.
\item There are \emph{uncountable} $\Sigma^1_2$ sets $A_0$ and $A_1$ such that $A_0\cup A_1=\R$ and for all $a\in \R$ the set $(a+A_0)\cap A_1$ is countable.
\item The plane can be covered by three $\Sigma^1_2$ clouds\footnote{see \S\ref{s.clouds} for the definition of clouds.} with centres in $L$.
\item There is a $\Sigma^1_2$ surjection $f:\R\to\R^2:x\mapsto (f_1(x),f_2(x))$ such that either $f_1'(x)$ or $f_2'(x)$ exists for all $x\in\R$.
\end{enumerate}
\end{theorem}

Here (2) and (3) correspond to CH equivalences proven by Sierpinski \cite{sierpinski65}; (4) to an equivalence due to Banach and Trzeciakiewicz, \cite{banach32,Trzeciakiewicz33}; (5) to an equivalence due to Komjath \cite{kope01}; and (6) to an equivalence proven by Morayne \cite{morayne84}.

The proofs of the above equivalences also offer an explanation for why and when a classical CH equivalence admits a $\Sigma^1_2$ counterpart. The reason that the above $\Sigma^1_2$ translations work can be found in the structure of the proofs of the corresponding classical CH equivalences. Though it is not always immediately clear from the literature, there is a common underlying structure of the proofs of CH from the given statement, and in fact of the statements themselves. Roughly speaking, the structure is as follows: The \emph{statements} are of the form that there exists certain sets (or $n$-ary relations) $R_1, R_2,\ldots$ and functions $f_1,f_2,\ldots$ which satisfy some finiteness or countability requirement, and that \emph{all} reals must satisfy some relations that are expressed in terms of the given sets and functions. The \emph{proof} that such a statement implies CH then can be cast in the following general form: One fixes a set of reals of size $\aleph_1$, and forms a ``hull'' of reals that satisfies the relevant relations with this fixed set of reals. The countability condition on the sets and functions $R_1,R_2,\ldots$, $f_1,f_2,\ldots$ then implies that this ``hull'' must have size $\aleph_1$. The statement is then seen to imply that in fact \emph{all} reals are in this hull, hence $2^{\aleph_0}=\aleph_1$.

In practice, one more often argues indirectly by assuming $\neg$CH, and then use this to produce a real which is ``transcendental'' in the sense that it fails to satisfy the prescribed relations. In the $\Sigma^1_2$ translations we consider, this corresponds to assuming that there is a non-constructible real. In our proofs the finiteness/countability conditions are then used, in conjunction with the Mansfield-Solovay perfect set theorem (see Theorem \ref{t.pst} below), to prove that the constructible reals are indeed a suitable ``hull''. Another important tool is the Shoenfield absoluteness theorem (see \cite{shoenfield61} or \cite[25.20]{jech03}), which allow us to work in a model of the form $L[x]$, $x\in\R$, which for the purpose of counting arguments can then be assumed to satisfy $\aleph_1^{L[x]}=\aleph_1^L$, see Lemma \ref{l.lav}.

Using the same ideas we also prove the following two partition relations which hold for $\Sigma^1_2$ colourings precisely when there are non-constructible reals.

\begin{theorem}\label{t.partition}
The following are equivalent:
\begin{enumerate}
\item $\R\nsubseteq L$.
\item For every $\Sigma^1_2$ function $f:\R\times\R\to\omega$ there are sets $C,D\subseteq\R$ such that $|C|=|D|=\aleph_0$ and $f\restrict C\times D$ is monochromatic.
\item For every $\Sigma^1_2$ colouring $g:\R\to\omega$ there are four distinct $x_{00},x_{01},x_{10},x_{11}\in\R$ of the same colour such that
$$
x_{00}+x_{11}=x_{01}+x_{10}.
$$
\end{enumerate}
\end{theorem}

This theorem, as well as Theorem \ref{t.mainthm1}, naturally relativizes to $L[a]$ and $\Sigma^1_2(a)$, where $a\in\R$ is a parameter.

\bigskip

The authors wish to thank Philip Welch for pointing out a mistake in a previous version of the paper, and for his comments regarding Remark \ref{r.endremark}.

\section{Definitions and preliminaries}

In this section we collect various general definitions and preliminary observations that are needed in our proofs. For this purpose, it is immensely practical to follow the (effective) descriptive set-theoretic convention and use $\R$ to stand for \emph{any} recursively presented uncountable Polish space (which is warranted since all such spaces are isomorphic by a $\Delta^1_1$ bijection, see \cite{moschovakis09}.) This convention will, however, cause problems later, where $\R$ will need to stand for the actual (linearly ordered field of) real numbers. Henceforth, we will use $\mathscr R$ to denote the descriptive set-theoretic \emph{reals} and $\R$ for the actual real line.

We shall assume that the reader is familiar with the basic elements of (effective) descriptive set theory, as found in e.g. \cite{manwei85}, \cite{moschovakis09} or \cite{drake74}, though we briefly review the most important notions below. Our notation is, for the most part, in line with that of \cite{moschovakis09}, and in particular, recursively presented Polish spaces are denoted with script letters $\mathscr X,\mathscr Y,\mathscr Z,\ldots$.

\subsection{$\Sigma^1_2$ sets and functions} In this paper, a $\Sigma^1_2$ set is a set that can be defined by a $\Sigma^1_2$ predicate, a $\Pi^1_2$ set is a set that is the complement of a $\Sigma^1_2$ set, and a $\Delta^1_2$ set is a set that is both $\Sigma^1_2$ and $\Pi^1_2$. We denote by $\Sigma^1_2(a)$, $\Pi^1_2(a)$ and $\Delta^1_2(a)$ the corresponding relativized pointclasses, where $a$ is some real (i.e., $a\in\mathscr R$.)

In this paper, we will say that a (total) function $f:\mathscr X\to\mathscr Y$ is $\Sigma^1_2$ (or, more generally, $\Sigma^1_2(a)$) if the graph of $f$ is a $\Sigma^1_2$ ($\Sigma^1_2(a)$) subset of $\mathscr X\times\mathscr Y$. If a function has a $\Sigma^1_2$ graph then in fact it is a $\Delta^1_2$ graph since if $\psi(x,y)$ is a $\Sigma^1_2$ predicate defining (the graph of) $f$ then
$$
f(x)=y\iff (\forall z)(\neg \psi(x,z)\vee z=y),
$$
which shows that $f$ has a $\Pi^1_2$ definition as well. We will say that a $\Sigma^1_2$ predicate $\psi(x,y)$ \emph{defines a function} if there is a total function $f:\mathscr X\to\mathscr Y$ such that $f(x)=y\iff\psi(x,y)$. The reader should be warned that this notion is sensitive to the model of set theory in which we work, since a predicate which defines a total function in one model may only define a partial function in another (for example, take $\psi(x,y)$ to be a $\Sigma^1_2$ predicate which says that $x=y$ and $x\in L$.) Note, however, that
\begin{equation}\label{eq.abs}
(\forall x)(\forall y,y') \psi(x,y)\wedge\psi(x,y')\implies y=y'
\end{equation}
is $\Pi^1_2$ and therefore absolute, and so if \eqref{eq.abs} is satisfied in one model, it is satisfied in all. In other words, a $\Sigma^1_2$ predicate which defines a partial function will do so in any model, but it may fail to define a total function in all models even if it does so in one.

\subsection{Coding the $L_\alpha$}\label{s.coding} Our notation follows that of \cite[p. 167ff.]{kanamori03}, with very
few differences. For convenience we recall the definitions and facts
that are most important for the present paper.

The canonical wellordering of $L$ will be denoted $<_L$. The
language of set theory (LOST) is denoted $\mathcal L_\epsilon$. If
$x\in 2^\omega$ then we define a binary relation on $\omega$ by
$$
m\inx{x}n\iff x(\langle m,n\rangle)=1,
$$
where $\langle \cdot,\cdot\rangle$ refers to some (fixed) standard G\"odel
pairing function of coding a pair of integers by a single integer.
We let
$$
M_x=(\omega,\epsilon_x)
$$
be the $\mathcal L_\epsilon$ structure coded by $x$. If $M_x$ is
wellfounded and extensional then we denote by $\tr(M_x)$ the
transitive collapse of $M_x$, and by $\pi_x:M_x\to\tr(M_x)$ the
corresponding isomorphism.

The following proposition encapsulates the basic descriptive
set-theoretic correspondences between $x$, $M_x$ and the
satisfaction relation. We refer to \cite[13.8]{kanamori03} and the
remarks immediately thereafter for a proof.

\begin{prop}
(a) If $\varphi(v_0,\ldots,v_{k-1})$ is a LOST formula with all free
variables shown then
$$
\{(x,n_0\ldots,
n_{k-1})\in 2^\omega\times\omega\times\cdots\times\omega: M_x\models
\varphi[n_0,\ldots,n_{k-1}]\}.
$$
is arithmetical.

(b) For $x\in 2^\omega$ such that $M_x$ is wellfounded and
extensional, the relation
$$
\{(m,y)\in\omega\times \mathscr R: \pi_x(m)=y\}
$$
is arithmetical in $x$.

(c) There is a LOST sentence $\sigma_0$ such that if
$M_x\models\sigma_0$ and $M_x$ is wellfounded and extensional, then
$M_x\simeq L_\delta$ for some limit ordinal $\delta<\omega_1$.

(d) There is a LOST formula $\varphi_0(v_0,v_1)$ which defines the
canonical wellordering $<_L$ of $L_\delta$ for all
$\delta>\omega$.\label{basicprop}
\end{prop}

\medskip

Define as in \cite[p.
170]{kanamori03} the restriction $M_x\restrict k$, for $x\in 2^\omega$
and $k\in\omega$, to be the $\mathcal L_{\epsilon}$ structure
$$
M_x\restrict k=(\{n: n\inx x k\},\epsilon_x).
$$

For $\mathscr X$ a recursively presented Polish space, let $R^{\mathscr X}\subseteq \mathscr X\times 2^\omega$ be defined by
\begin{align*}
R^{\mathscr X}(y,x)\iff M_x\text{ is well-founded, extensional, and } M_x\models\sigma_0\wedge (\exists n) y=\pi_x(n)\\
(\forall z<_Lx) (M_z\text{ is well-founded, extensional and } M_z\models\sigma_0)\implies (\forall k) \pi_z(k)\neq y
\end{align*}
In other words, $R^{\mathscr X}(y,x)$ holds iff $x$ is the least code for an $L_\alpha$, $\alpha$ a limit, such that $y\in L_\alpha$. The relation $R^{\mathscr X}$ is $\Delta^1_2$.

\subsection{Coding initial segments.}\label{s.initseg} Let $\prec$ denote $<_L\restrict \mathscr R$, the canonical well-ordering of $\mathscr R$ in $L$. This is a $\Delta^1_2$ wellordering which has a good coding of initial segments. More precisely, $\prec$ is a \emph{strongly} $\Delta^1_2$ well-ordering, which means that $\prec$ has length $\omega_1$ and $\IS\subseteq\mathscr R\times\mathscr R^{\leq\omega}$ defined by
$$
\IS(x,v)\iff (\forall z\prec x)(\exists n) v(n)=z\wedge (\forall i,j) i=j\vee v(i)\neq v(j)
$$
is $\Delta^1_2$. The point is that quantifications over an initial segments of $\prec$  can be replaced by a quantifier over $\omega$ in hierarchy calculations, see \cite[5A.1]{moschovakis09} for details. We also define a function $\IS^*:\mathscr R\to\mathscr R^{\leq\omega}$ and a partial function $\IS^\#:\mathscr R\times\mathscr R\to\omega$ by
\begin{align*}
&\IS^*(x)=v\iff \IS(x,v)\wedge (\forall w\prec v) \neg \IS(x,w)\\
&\IS^\#(x,y)=n\iff \IS^*(x)(n)=y.
\end{align*}
These are $\Sigma^1_2$.

\subsection{The size of $L\cap\R$.} There are several counting arguments below that rely on having some information about the cardinality of sets of reals in $L$. The following simple observations is extremely useful for this purpose:

\begin{lemma}\label{l.lav}
(1) If there is a non-constructible real in $V$, then there is a non-constructible real $x\in V$ such that $\aleph_1^{L[x]}=\aleph_1^L$.

(2) Suppose $\psi$ is a $\Sigma^1_2(a)$ predicate defining the set $A$, where $a\in L$. Then if $A$ is uncountable, then $A\cap L$ is uncountable in $L$.
\end{lemma}

\begin{proof}
(1) If $\aleph_1^V=\aleph_1^L$, then any non-constructible $x\in V$ will do. If $\aleph_1^L$ is countable in $V$, then there must be a real $x\in V$ which is Cohen over $L$. For any such $x$ it holds that $\aleph_1^{L[x]}=\aleph_1^L$.

(2) If $A\cap L$ is countable in $L$ then there is some $v:\N\to\mathscr R$ in $L$ such that
$$
(\forall x) (\psi(x)\longrightarrow (\exists n) v(n)=x)
$$
holds. Since this is $\Pi^1_2(a,v)$ it is absolute, and so $A$ is countable.
\end{proof}

The typical application of (1) above will be that if we know that some statement which is downwards absolute holds in $V$, and $\R\nsubseteq L$, then the statement holds in some $L[x]$ where $x\notin L$, and the constructible reals have cardinality $\aleph_1$ in $L[x]$.

Finally, we recall the perfect set theorem for $\mathbf{\Sigma}^1_2$ sets by Mansfield and Solovay which will be used often:

\begin{theorem}[Mansfield \cite{mansfield70}, Solovay \cite{solovay69}]\label{t.pst}
Let $A$ be a $\Sigma^1_2(a)$ set. Then either $A\subseteq L[a]$, or else $A$ contains a perfect set. In particular, if a $\Sigma^1_2$ set contains a non-constructible real then it is uncountable.
\end{theorem}

\section{Results}

\subsection{Sierpinski's equivalences}

In this section we consider the $\Sigma^1_2$ counterparts of two of Sierpinski's classical CH equivalences (see e.g. \cite{sierpinski65}). The first is the counterpart to: CH is equivalent to the existence of two sets $A,B\subseteq \R^2$ with $A\cup B=\R^2$ such that all vertical sections of $A$ are countable and all horizontal sections of $B$ are countable.

We include a version of this that is stated in terms of covering the plane by graphs of countably many functions, since this is needed later in section \ref{ss.morayne} below.

\begin{theorem}\label{t.sierpinski1}
The following are equivalent:

\begin{enumerate}
\item $\R\subseteq L$.
\item There is a $\Sigma^1_2$ linear order $<$ of $\R$ such that for all $x\in\R$ the initial segment $\{y\in\R: y<x\}$ is countable.
\item There are $\Sigma^1_2$ sets $A,B\subseteq\R^2$ such that $A\cup B=\R^2$, and all the sections $A_x=\{y\in \R: (x,y)\in A\}$ and $B^y=\{x\in \R: (x,y)\in B\}$ are countable.
\item There are $\Sigma^1_2$ functions $F_A:\R\times\omega\to\R$ and $F_B:\R\times\omega\to\R$ such that $A=\{(x,F_A(x,n)): x\in\R,n\in\omega\}$ and $B=\{(F_B(y,n),y):y\in\R,n\in\omega\}$ satisfy (3).

\end{enumerate}
\end{theorem}

\begin{proof}

$(1)\implies (4)$. Let $z$ be the $\prec$-least element with an infinite initial segment. Let 
$$
F_A(x,n)=\IS^*({\max}_{\prec}(x,z))(n)
$$
and 
$$
F_B(x,n)=\left\{\begin{array}{ll}
x & \text{if } n=0\\
\IS^*(\max_{\prec}(x,z))(n-1) & \text{if } n>0
\end{array}\right.
$$
where $\max_{\prec}(x,z)$ is the larger of $x$ and $z$ in $\prec$.

$(4)\implies (3)$ is clear. For $(3)\implies (1)$, suppose that there is $x_0\in\R\setminus L$ but that (3) holds. By Lemma \ref{l.lav} we may assume that $\aleph_1^{L[x_0]}=\aleph_1^L$ and that $V=L[x_0]$, since if $(3)$ holds it holds in $L[x_0]$. Since the section $A_{x_0}$ is countable we can find $y\in (\R\cap L)\setminus A_{x_0}$, and so $(x_0,y)\in B$ since $A\cup B=\R^2$. But this means that $B^y$, which is a $\Sigma^1_2(y)$ set, contains a non-constructible real (namely $x_0$), and so since $y\in L$ it follows by the perfect set theorem (Theorem \ref{t.pst}) that it must be uncountable, a contradiction.

Finally, $(1)\implies(2)$ is clear, since the canonical well-ordering of $\R$ in $L$ satisfies (2), and $(2)\implies(3)$ follows since defining $A=\{(x,y)\in \R^2: y<x\}$  and $B=\{(x,y): x\leq  y\}$ clearly works. 
\end{proof}

Next we consider the $\Sigma^1_2$ counterpart to the following CH equivalence due to Sierpinski (see \cite{sierpinski65}): CH holds iff there are sets $A_1,A_2,A_3\subseteq\R^3$ such that $A_1\cup A_2\cup A_3=\R^3$, and every line $l$ in the direction of the $x_i$-axis meets $A_i$ in finitely many points.

\begin{theorem}
All reals are constructible if and only if there are $\Sigma^1_2$ sets $A_1,A_2,A_3\subseteq\R^3$ such that $A_1\cup A_2\cup A_3=\R^3$, and every line $l$ in the direction of the $x_i$-axis meets $A_i$ in finitely many points.
\end{theorem}

\begin{proof}
Suppose all reals are constructible. Define, for $i=1,2,3$, the set $\tilde A_i$ by
\begin{align*}
(x_1,x_2,x_3)\in \tilde A_i \iff  &\text{if $x_j={\max}_\prec\{x_1,x_2,x_3\}$ then $x_j\neq x_i$,}\\
&\text{and if $k\neq i,j$ then $\IS^\#(x_j,x_i)<\IS^\#(x_j,x_k)$}\\
\iff & (\forall j) (x_j={\max}_\prec\{x_1,x_2,x_3\}\implies (x_j\neq x_i\wedge\\
&((\forall k\leq 3)(k\neq i\wedge k\neq j)\implies \IS^\#(x_j,x_i)<\IS^\#(x_j,x_k))
\end{align*}
Clearly $\tilde A_i$ is $\Delta^1_2$. Let $A_i=\tilde A_i\cup \{(x,y,z):x=y=z\}$. Then $\R^3=A_1\cup A_2\cup A_3$. If $l$ is a line parallel to an axis, say $l=\{(x,b,c):x\in\R\}$, then by definition there are only finitely many $x$ such that $(x,b,c)\in A_1$.

For the converse, suppose there is $x_0\in \R\setminus L$. As before, we may assume that $V=L[x_0]$ and that $\aleph_1^{L[x_0]}=\aleph_1^L$. If $(u,v)\in \R^2\cap L$, then the line  $\{(u,v,x):x\in\R\}\cap A_3$ is a finite $\Sigma^1_2(u,v)$ set, and so by Theorem \ref{t.pst} it does not contain a non-constructible real. Thus $(u,v,x_0)\notin A_3$ for all $u,v\in\R\cap L$. For any $u\in\Q$ the set $\{(u,x,x_0):x\in\R\}\cap A_2$ is finite, and so since $\aleph_1^L=\aleph_1$ there must be some $x_1\in\R\cap L$ such that $(u,x_1,x_0)\notin A_2$ for all $u\in\Q$. Since $A_1\cap\{(u,x_1,x_0):u\in\R\}$ is finite, it follows that there is $x_2\in\Q$ such that $(x_2,x_1,x_0)\notin A_1\cup A_2\cup A_3$.
\end{proof}

\subsection{Banach-Trzeciakiewicz's equivalence}

\cite{banach32} and \cite{Trzeciakiewicz33} contain the following equivalence: CH holds if and only if there are uncountable sets $A_0,A_1\subseteq\R$ such that $A_0\cup A_1=\R$ and for each $a\in\R$ the set $(a+A_0)\cap A_1$ is countable. We have the following $\Sigma^1_2$ counterpart:

\begin{theorem}
All reals are constructible if and only if there are \emph{uncountable} $\Sigma^1_2$ sets $A_0$ and $A_1$ such that $A_0\cup A_1=\R$ and for all $a\in \R$ the set $(a+A_0)\cap A_1$ is countable.
\end{theorem}

\begin{proof}
If $\R\subseteq L$, it is easy to see that there is a $\Delta^1_2$ Hamel basis $H\subseteq\R$ for $\R$. (In fact, by \cite{miller89} there even is a $\Pi^1_1$ Hamel basis for $\R$.) Define a function $f:\R\to \R^{<\omega}$ by
$$
f(x)=(x_1,\ldots, x_n)\iff x_1,\ldots, x_n\in H\wedge x_1\prec\cdots\prec x_n\wedge (\exists q_1,\ldots, q_n\in\Q\setminus\{0\}) x=\sum_{i=1}^n q_ix_i.
$$
Clearly $f$ is $\Delta^1_2$. Write $H=H_0\cup H_1$, where $H_0$ and $H_1$ are disjoint uncountable $\Delta^1_2$ sets, and define
$$
x\in A_i\iff (\exists (x_1,\ldots, x_n)\in \R^{<\omega}) f(x)=(x_1,\ldots,x_n)\wedge x_n\in H_i.
$$
Then $A_i$ is $\Sigma^1_2$ (in fact, $A_i$ is $\Delta^1_2$) and $A_0\cup A_1=\R$. Fix $a\in\R$, and note that if $\max f(a)\prec \max f(x)$ and $x\in A_0$ then $a+x\in A_0$. Thus $(a+A_0)\cap A_1$ is countable since $\{x\in\R: \max f(x)\preceq \max f(a)\}$ is.

For the converse, suppose that there is a non-constructible real $x_0\in\R\setminus L$. By Lemma \ref{l.lav}.(1) we may assume that $\aleph_1^{L[x_0]}=\aleph_1^L$. From this and Lemma \ref{l.lav}.(2) it follows that $A_0\cap L$ and $A_1\cap L$ are uncountable, and so $A_0\cap L[x_0]$ and $A_1\cap L[x_0]$ are uncountable in $L[x_0]$. By assumption, for each $a\in A_0$ we either have $a+x_0\in A_0$ or $a+x_0\in A_1$. If the latter held for uncountably many $a\in A_0\cap L$ then $(x_0+A)\cap A_1$ would be uncountable, contrary to our assumption. Thus we can find $a\in L\cap A_0$ such that $a+x_0\in A_0$. Similarly, there is $b\in A_1\cap L$ such that $b+x_0\in A_1$. But since $b+x_0=a+x_0+(b-a)$ we now have that $b+x_0\in ((b-a)+A_0)\cap A_1$, and so this set, which is $\Sigma^1_2(b-a)$, contains the non-constructible real $b+x_0$, and so is uncountable.
\end{proof}

\subsection{Komjath's clouds}\label{s.clouds}

A \emph{cloud} in $\R^2$ is a set $A\subseteq\R^2$ such that for some point $\vec x\in\R^2$ (called a \emph{centre} of $A$) it holds that each infinite ray from $\vec x$ meets $A$ in at most finitely many points. In \cite{kope01} the following was shown:

\medskip

\noindent {\bf Theorem.} (Komjath). {\it CH is equivalent to that the plane can be covered by three clouds.}

\begin{theorem}
$\R\subseteq L$ is equivalent to that the plane can be covered by three $\Sigma^1_2$ clouds with centres in $L$.
\end{theorem}

\begin{proof}
Assume that $\R\subseteq L$. We will give $\Sigma^1_2$ definitions of clouds $A_0$, $A_1$ and $A_2$ centered at $a_0=(0,1)$, $a_1=(1,0)$ and $a_2=(0,0)$, respectively, such that $\R=A_0\cup A_1\cup A_2$. For $y\in\R^2\setminus\{a_0,a_1,a_2\}$ let $\overline{a_iy}$ denote the infinite ray starting at $a_i$ extending through $y$. Let $\mathscr E$ be the set of all infinite rays from $a_0$, $a_1$ or $a_2$. The set of $\mathscr E$ can be identified with the union of the three disjoint circles centered at $a_0,a_1$ and $a_2$, and so $\mathscr E$ is a recursively presented Polish space in a natural way. Let $\mathscr E_\alpha=\mathscr E\cap L_{\alpha}$. 

We define the set $A_i'\subseteq \R^2\times 2^\omega$ as follows: $(y,x)\in A_i'$ if and only if
\begin{enumerate}
\item $R^{\mathscr E}(\overline{a_iy},x)$, i.e., $x$ is $\prec$-least such that $M_x\simeq L_\alpha$ for the smallest limit $\alpha>\omega$ such that $\overline{a_iy}\in L_\alpha$.
\item If $(j_l)_{l\in\omega}$ is a strictly increasing sequence enumerating the set
$$
\{j\in\omega:\pi_x(j)\in \mathscr E_\alpha\setminus \bigcup \{\mathscr E_\delta:\delta<\alpha,\delta\text{ a limit}\}\}
$$
and the ray $\overline{a_iy}$ is $\pi_x(j_l)$, then $y$ is a point of intersection between $\pi_x(j_l)$ and one of the rays $\pi_x(j_0),\ldots,\pi_x(j_{l-1})$ or $\overline{a_ja_k}$, $j\neq k$ and $j,k\neq i$.
\end{enumerate}
Then $A_i'$ is $\Sigma^1_2$ since (2) can (given that (1) holds) be expressed by saying (where $j,k\neq i$)
$$
(\exists l)[ \pi_x(l)=\overline{a_iy}\wedge((y\in \overline{a_iy}\cap \overline{a_{j_0}a_{j_1}})\vee ((\exists j< k) R^{\mathscr E}(\pi_x(j),x)\wedge y\in\pi_x(j)\cap\overline{a_iy}))].
$$
Let $A_i=\{y\in \R^2: (\exists x) A_i'(y,x)\}$, which clearly is a $\Sigma^1_2$ set, and note that if $y\in\R^2\cap L$ then there must be some $i\in\{0,1,2\}$ such that $y\in A_i$, and so $A_0\cup A_1\cup A_2=\R^2\cap L$, as required.

For the converse, assume that there are $\Sigma^1_2$ clouds $A_0, A_1$ and $A_2$ with centres in $L$ covering the plane. After possibly applying an affine transformation (defined in $L$), we may assume that $A_0, A_1$ and $A_2$ are centered at $(0,1), (1,0)$ and $(0,0)$, respectively.

By the usual arguments, we can assume that $V=L[r]$ for some $r\in \R\setminus L$ and that $\aleph_1=\aleph_1^L$. Define an equivalence relation in $(0,\frac \pi 4)$ by 
$$
\alpha\sim\alpha'\iff \frac {1-\tan(\alpha')}{1-\tan(\alpha)}\in\Q_+.
$$
Then $\sim$ has countable classes and $\alpha\in L$ iff $[\alpha]_\sim\subseteq L$.

For $\alpha,\beta\in (0,\frac\pi 4)$, let $l_\alpha$ denote the straight line in the plane given by the equation $\tan(\alpha)x+y=1$, and $t_\beta$ be the line given by $x+\tan(\beta)y=1$. Note that the intersection point $(x,y)$ of $l_\alpha$ and $t_\beta$ satisfies $\frac y x=\frac{1-\tan(\alpha)}{1-\tan(\beta)}$.

Consider $\alpha\in (0,\frac\pi 4)\cap L$. Since $l_\alpha\cap A_0$ is a finite $\Sigma^1_2(\alpha)$ set, it cannot contain any non-constructible points (by Theorem \ref{t.pst}, for example.) Thus if $\beta\in (0,\frac\pi 4)\setminus L$, then the intersection of $l_\alpha$ and $t_\beta$ cannot be in $A_0$. So fix $\beta_0\in (0,\frac\pi 4)\setminus L$. Since $A_1\cap\{t_\beta:\beta\in [\beta_0]_\sim\}$ is countable there must be some $\alpha_0\in (0,\frac\pi 4)\cap L$ such that $t_\beta\cap t_\alpha\nsubseteq A_1$
for all $\beta\in [\beta_0]_\sim$ and $\alpha\in[\alpha_0]_\sim$, whence $t_\beta\cap t_\alpha\subseteq A_2$ for such $\alpha$ and $\beta$. For $n\in\N$, choose $\alpha_n\in [\alpha_0]_\sim$ and $\beta_n\in[\beta_0]_\sim$ such that
$$
\frac {1-\tan(\alpha_0)}{1-\tan(\alpha_n)}=n=\frac {1-\tan(\beta_0)}{1-\tan(\beta_n)}.
$$
Then for all $n\in\omega$ the intersection point $(x_n,y_n)\in l_{\alpha_n}\cap t_{\beta_n}$ satisfies
$$
\frac{y_n}{x_n}=\frac{1-\tan(\alpha_n)}{1-\tan(\beta_n)}=\frac{1-\tan(\alpha_0)}{1-\tan(\beta_0)}.
$$
and so they are all on the same line through $(0,0)$, and since $(x_n,y_n)\in A_2$ for all $n\in\N$ this contradicts that each ray from $(0,0)$ meets $A_2$ in finitely many points.
\end{proof}

\begin{remark}
It is interesting to note that in the previous proof, the assumption that $A_1$ and $A_2$ are $\Sigma^1_2$ were never used. Thus we have:
\end{remark}

\begin{corollary}
$\R\subseteq L$ is equivalent to that the plane can be covered by three clouds with centres in $L$, one of which is $\Sigma^1_2$.
\end{corollary}

\subsection{Differentiable functions after Morayne}\label{ss.morayne}

A \emph{Peano function} is a surjection $f:\R\to\R\times\R$. In \cite{morayne84}, Morayne proved that CH is equivalent to the existence of a Peano function $f(x)=(f_1(x),f_2(x))$ such that at every $x\in\R$ at least one of the derivatives $f_1'(x)$ or $f_2'(x)$ exists. We obtain the following corresponding $\Sigma^1_2$ version:

\begin{theorem}\label{t.diffsurj}
The following are equivalent:

\begin{enumerate}

\item All reals are constructible

\item There is a $\Sigma^1_2$ surjection $f:\R\to\R^2:x\mapsto (f_1(x),f_2(x))$ such that either $f_1'(x)$ or $f_2'(x)$ exists for all $x\in\R$.

\end{enumerate}
\end{theorem}

\begin{proof}[Proof of (1)$\implies$ (2)]
We will show that the construction from CH due to Morayne translates to the $\Sigma^1_2$ setting.

For this, first define $f_1(t)=t\sin(t)$ on $t\in (-\infty,1)=I_1$ and $f_2(t)=t\sin(t)$ on $t\in (-1,\infty)=I_2$. The sets
$$
C^i=\{(r,t)\in\R\times I_i:f_i(t)=r\},
$$
$i=1,2$, are $\Delta^1_1$ and for each $r\in\R$ the section $C^i_r=\{t\in I_i: (r,t)\in C^i\}$ is countably infinite. It follows from (the effective version of) the Lusin-Novikov Theorem \cite[18.10]{kechris95} that there are $\Delta^1_1$ functions $g_i:\R\to I_i^\omega$ such that $g_i(r)$ enumerates $C^i_r$ injectively. Now let $F_A$ and $F_B$ be the functions from Theorem \ref{t.sierpinski1}.(4), and define for $t\in\R\setminus I_1$
$$
f_1(t)=y\iff (\exists r\in\R)(\exists n\in\omega) f_2(t)=r\wedge g_2(r)(n)=t\wedge F_B(r,n)=y
$$
and for $t\in\R\setminus I_2$
$$
f_2(t)=y\iff (\exists r\in\R)(\exists n\in\omega) f_1(t)=r\wedge g_1(r)(n)=t\wedge F_A(r,n)=y.
$$
Note that whenever $t\notin I_1$ and $f_2(t)$ assumes the value $r$ for the $n$'th time as enumerated by $g_2(r)$, then $(f_1(t),f_2(t))=(F_B(r,n), r)$, and so the graph of $t\mapsto (f_1(t),f_2(t))$ covers $B=\{(F_B(r,n),r):r\in\R,n\in\omega\}$ as $t$ ranges in $I_1$. Similarly, the graph of $t\mapsto (f_1(t),f_2(t))$ covers $A=\{(r,F_A(r,n)):r\in\R, n\in\omega\}$ as $t$ ranges over $I_2$. Thus $t\mapsto (f_1(t), f_2(t))$ is a $\Sigma^1_2$ Peano function with $f_1$ differentiable on $I_1$ and $f_2$ differentiable on $I_2$.
\end{proof}

The proof of $(2)\implies (1)$ in Theorem \ref{t.diffsurj} requires several lemmata. We start with a general observation about open $\Pi^1_2$ sets. Recall that the class of $\Pi^1_2$ sets is \emph{$\omega$-parametrized}, meaning that for any recursively presented Polish $\mathscr X$, there is a $\Pi^1_2$ set $P^{(\mathscr X)}\subseteq\omega\times\mathscr X$ such that 
$$
P_n^{(\mathscr X)}=\{x\in\mathscr X:(n,x)\in P\}
$$
enumerates the $\Pi^1_2$ sets in $\mathscr X$. In particular, there is such a set $P^{(\omega)}\subseteq\omega\times\omega$ parametrizing the $\Pi^1_2$ subsets of $\omega$. We let 
$$
\mathfrak a=\{\langle n,m\rangle: (n,m)\in P^{(\omega)}\},
$$
where $\langle\cdot,\cdot\rangle$ is some standard G\"odel pairing function. Note that $\mathfrak a\in L$.

\begin{lemma}\label{l.openpi-1-2}
Suppose $A\subseteq\mathscr X$ is an \emph{open} $\Pi^1_2$ set. Then there is a $\Sigma^0_1(\mathfrak a)$ predicate $\psi(x)$ such that $x\in A\iff \psi(x)$.
\end{lemma}
\begin{proof}
Let $d$ be a compatible metric on $\mathscr X$ and let $(x_n)_{n\in\omega}$ be a dense sequence in $\mathscr X$ such that $(d,(x_n)_{n\in\omega})$ is a recursive presentation of $\mathscr X$. Let $(q_m)_{m\in\omega}$ enumerate (effectively) the positive rationals, and define
$$
a=\{\langle n,m\rangle\in\omega: (\forall x) d(x,x_n)<q_m\implies x\in A\}.
$$
Then the set $a\subseteq\omega$ is $\Pi^1_2$, and
$$
x\in A\iff (\exists n,m) \langle n,m\rangle\in a\wedge d(x,x_n)<q_m
$$
which gives a $\Sigma^0_1(a)$ definition of $A$, whence $A$ is $\Sigma^0_1(\mathfrak a)$.
\end{proof}

\begin{lemma}\label{l.defosc}
Let $\psi(x,y)$ be a $\Delta^1_2$ predicate which defines a function $f: \R\to\R$. Then:

(1) There is a $\Pi^0_2(\mathfrak a)$ predicate $\phi(x)$ such that in any model in which $\psi$ defines a function we have: $\phi(x)$ holds if and only if $x$ is a point of continuity of $f$.

(2) There is a $\Pi^0_2(b)$ predicate $\hat\psi(x,y)$ with parameter $b\in L$ such that in any model where $\psi$ defines a function we have: $\hat\psi(x,y)$ if and only if $\psi(x,y)\wedge\phi(x)$.
\end{lemma}

\begin{proof}
(1) Recall that for $x\in\R$, the \emph{oscillation} of $f$ at $x$ is defined as
$$
\osc_f(x)=\inf\{\diam(f(U)): x\in U\wedge U\subseteq\R\text{ is open}\},
$$
and that $x$ is a point of continuity precisely when $\osc_f(x)=0$. Let $\phi(x,\varepsilon)$ be the following predicate:
\begin{align*}
&(\exists q,r,\delta\in\Q_+) |x-q|<r\wedge [(\forall x_0,x_1)(\forall y_0,y_1) (f(x_0)=y_0\wedge f(x_1)=y_1\wedge\\
& |x_0-q|<r\wedge |x_1-q|<r)\longrightarrow |y_0-y_1|<\varepsilon-\delta].
\end{align*}
This is $\Pi^1_2$ and $\phi(x,\varepsilon)$ holds precisely when $\osc_f(x)<\varepsilon$. On the other hand, it is easy to see that $\{(x,\varepsilon)\in\R\times\Q_+:\osc_f(x)<\varepsilon\}$ is open (when $\Q_+$ has the discrete topology), and so $\{(x,\varepsilon)\in\R\times\Q_+:\hat\psi(x,\varepsilon\}$ is an open $\Pi^1_2$ set. It follows from Lemma \ref{l.openpi-1-2} that there is a $\Sigma^0_1(\mathfrak a)$ predicate $\hat\phi(x,\varepsilon)$ such that $\hat\psi(x,\varepsilon)$ iff $\hat\phi(x,\varepsilon)$. Thus if we let $\phi(x)$ be $(\forall\varepsilon\in\Q_+) \hat\phi(x,\varepsilon)$ then $\phi(x)$ is a $\Pi^0_2(\mathfrak a)$ predicate which holds precisely when $x$ is a point of continuity of $f$, and $\phi$ does so in any model where $\psi(x,y)$ defines a function.

(2) Fix a sequence $(x_n)_{n\in\omega}$ in $\R\cap L$ such that 
$$
L\models \text{``$(x_n)_{n\in\omega}$ is dense in $\{x\in\R:\phi(x)\}$''}.
$$
To say that $(x_n)$ is dense in $\{x\in\R:\phi(x)\}$ can be expressed as
$$
(\forall\varepsilon\in\Q_+)(\forall x) (\phi(x)\longrightarrow (\exists n) |x_n-x|<\varepsilon),
$$
which is $\Pi^1_1(\mathfrak a,(x_n)_{n\in\omega})$, and so this statement is absolute. Let $(y_n)_{n\in\omega}$ be the sequence in $\R\cap L$ defined by $y_n=f(x_n)$, and let $\hat\psi(x,y)$ be the predicate
$$
\phi(x)\wedge (\forall\varepsilon\in\Q_+) (\exists n) |x_n-x|<\varepsilon\wedge |y_n-y|<\varepsilon.
$$
Then $\hat\psi(x,y)$ is $\Pi^0_2(\mathfrak a,(x_n)_{n\in\omega},(y_n)_{n\in\omega})$ and since $f$ is continuous on the set $\{x\in\R:\phi(x)\}$ it holds that
$$
\hat\psi(x,y)\iff \phi(x)\wedge \psi(x,y)
$$
in any model where $\psi$ defines a function, as required.
\end{proof}

\begin{lemma}\label{l.setH}
Let $f:\R\to\R$ be a function. Then:

(1) There is a $\mathbf{\Pi}^1_1$ set $H\subset\R$ such that
$$
\{x\in\R: f'(x)\text{ exists}\}\subseteq H
$$
and $\{y\in\R: f^{-1}(y)\cap H\text{ is uncountable}\}$ is Lebesgue null.

(2) If $f$ is defined by the $\Delta^1_2$ predicate $\psi(x,y)$ then there is a $\Pi^1_1(b)$ predicate $\chi(x)$ with a parameter $b\in L$ such if we let $H=\{x\in\R:\chi(x)\}$ then (1) holds for this $H$ and $f$ defined by $\psi$ in any model where $\psi(x,y)$ defines a function.

\end{lemma}
\begin{proof}
(1) Let $C$ be the set of points of continuity of $f$. It is well-known that this is a $G_\delta$ set. Define
\begin{equation}\label{eq.def}
x\in H\iff x\in C\wedge (\exists y)(\forall\varepsilon>0)(\exists\delta>0)(\forall z\in C\setminus\{x\})\left[|x-z|<\delta\implies |\frac{f(x)-f(z)}{x-z}-y|<\varepsilon\right].
\end{equation}
It is clear that if $f'(x)$ exists then $x\in H$.

\begin{claim}
$H$ is $\mathbf{\Pi}^1_1$.
\end{claim}
\begin{proof}
Let $\bar f:\R\to\R$ be a Borel function such that $\bar f\restrict C=f\restrict C$. We claim that $x\in H$ if and only if
\begin{equation}\label{eq.pi11def}
x\in C\wedge (\forall\varepsilon>0)(\exists q\in\Q)(\exists\delta>0)(\forall z\neq x)\left[(z\in C\wedge|x-z|<\delta)\implies\left|\frac{\bar f(x)-\bar f(z)}{x-z}-q\right|<\varepsilon\right].
\end{equation}
If $x$ is isolated in $C$ then clearly \eqref{eq.def} holds for $x$ if and only if \eqref{eq.pi11def} holds. So assume that $x$ is not isolated. If \eqref{eq.pi11def} holds for $x$, let $q_n$ witness \eqref{eq.pi11def} with $\varepsilon=\frac 1 {2^{n+1}}$. Then $|q_{n+1}-q_n|\leq\frac 1 {2^n}$ so $q_n$ is Cauchy, and if we let $y=\lim_{n\to\infty} q_n$ then $y$ is easily seen to be a witness to \eqref{eq.def}. Conversely, if \eqref{eq.def} holds for $x$, and $y$ is a witness to this, then let $q_n\in\Q$ be a sequence of rationals such that $q_n\to y$. Then it is clear that for all $\varepsilon>0$ we can find some $n$ such that $q_n$ is a witness to that \eqref{eq.pi11def} holds.
Since \eqref{eq.pi11def} is $\mathbf{\Pi}^1_1$, the claim is proved.
\end{proof}

\begin{claim}
$\{y\in\R: f^{-1}(y)\cap H\text{ is uncountable}\}$ is Lebesgue null.
\end{claim}
\begin{proof}
The proof uses the idea from \cite[Ch. 5.15]{koto06}. It clearly suffices to show for all $m\in\N$ that the sets
$$
Y_m=\{y\in [-m,m]: f^{-1}(y)\cap H\text{ is uncountable}\}
$$
are null; we will prove this for $Y_1$, from which the other cases follow by rescaling the codomain of $f$ (or by an identical proof.) For $y\in Y_1$, pick $t_y\in f^{-1}(y)\cap H$ such that $f'(t_y)=0$. Such a $t_y$ exists since when $y\in Y_1$ the set $f^{-1}(y)\cap H$ is uncountable and so it contains an accumulation point, and as $f$ is constant on this set we must have $f'(t)=0$ at any accumulation point. Let $T=\{t_y:y\in Y_1\}$, and note that $f(T)=Y_1$.

Let $\varepsilon>0$, and for each $t\in T$ let $1>\delta_t>0$ be such that for all $z\in C$ with $|t-z|<\delta_t$ we have
$$
\left|\frac{f(t)-f(z)}{t-z}\right|<\varepsilon,
$$
and let $I_t=(t-\delta_t,t+\delta_t)$. Note that for any $z\in I_t\cap C$ we have $f(z)\in (f(t)-\varepsilon\delta_t,f(t)+\varepsilon\delta_t)$, and so we have $\mu(f(I_t\cap C))\leq 2\varepsilon\delta_t$. Since the intervals $I_t$ cover $T$, we can find $t_i\in T$, $i\in\N$, such that $U=\bigcup_{t\in T} I_t=\bigcup_{i\in\N} I_{t_i}$. We claim that $\mu(f(U\cap C))\leq 4\varepsilon$. To see this it is enough to prove that $\mu(f(K\cap C))\leq4\varepsilon$ for all compact $K\subseteq U$. If $K\subseteq U$ is compact, then we can find $N\in\N$ such that $K\subseteq\bigcup_{i=1}^N I_{t_i}$. Moreover, after possibly going to a subcover, we can assume that each $x\in K$ is contained in at most two different intervals $I_{t_i}$, and so we have $\sum_{i=1}^N 2\delta_{t_i}=\sum_{i=1}^N\mu(I_{t_i})\leq 4$. Thus
$$
\mu(f(C\cap K))\leq\mu(f(\bigcup_{i=1}^N I_{t_i}\cap C))\leq\sum_{i=1}^N\mu(f(I_{t_i}\cap C))\leq \sum_{i=1}^N 2\varepsilon\delta_{t_i}=4\varepsilon,
$$
as required. It follows that $\mu(f(U\cap C))=0$, and so since $T\subseteq U\cap C$ we have $
\mu(Y_1)=\mu(f(T))=0$
\end{proof}

(2) Let $\phi(x)$ and $\hat\psi(x,y)$ be the predicates defined in Lemma \ref{l.defosc}, and let $\chi(x)$ be the predicate
\begin{align*}\label{eq.pi11adef}
&\phi(x)\wedge(\forall\varepsilon>0)(\exists q\in\Q)(\exists\delta>0)(\forall z)(\forall y_0,y_1)\\
&(\phi(z)\wedge z\neq x\wedge|x-z|<\delta\wedge\hat\psi(x,y_0)\wedge\hat\psi(z,y_1))\longrightarrow |\frac{y_0-y_1}{x-z}-q|<\varepsilon.
\end{align*}
Then $\chi(x)$ is $\Pi^1_1(b)$ (where $b\in L$ is the parameter in $\hat\psi$), and if $\psi(x,y)$ defines a function then the set $\{x\in\R:\chi(x)\}$ is equal to the set $H$ defined in \eqref{eq.def}.
\end{proof}

\begin{proof}[Proof of Theorem \ref{t.diffsurj}]
We may assume that $\aleph_1^L=\aleph_1$. Fix $f:\R\to\R^2:x\mapsto (f_1(x),f_2(x))$ as in the statement of the theorem. Applying Lemma \ref{l.setH} to $f_1$ and $f_2$, there are $\Pi^1_1(b)$ ($b\in L$) sets $H_1$ and $H_2$ defined by $\Pi^1_1(b)$ formulas $\chi_1(x)$ and $\chi_2(x)$ such that 
$$
Y_i=\{y\in\R: |f_i^{-1}(y)\cap H_i|>\aleph_0\}
$$
is Lebesgue null, and such that the points of differentiability of $f_i$ are contained in $H_i$. Let
$$
Y_i^L=\{y\in\R\cap L: L\models |f_i^{-1}(y)\cap H_i|>\aleph_0\}.
$$
We claim that $Y_i\cap L=Y_i^L$. Since $\aleph_1=\aleph_1^L$ it is clear that if $y\in Y_i^L$ then $y\in Y_i$. On the other hand, note that the set $\Gamma_{i,y}=f^{-1}_i(y)\cap H_i$ is $\Delta^1_2(b,y)$, and so if $y\in L$ then by Lemma \ref{l.lav} the set $\Gamma_{i,y}$ is countable if $\Gamma_{i,y}\cap L$ is. So if $y\in L\setminus Y_i^L$ then $y\notin Y_i$, as required.

Let $\R^*=\R\setminus(Y_1\cup Y_2)$, which has full measure, and let $A_1=f(H_2)$ and $A_2=f(H_1)$. Then $A_1$ and $A_2$ are $\Sigma^1_2(b)$ sets, and since either $f'_1(t)$ or $f'_2(t)$ exists for all $t\in\R$ we must have that $\R=A_1\cup A_2$. For any $r\in\R^*$ the set $f_i^{-1}(r)\cap H_i$ is countable by the definition of $\R^*$, and so there are at most countably many $t\in H_i$ such that $f_i(t)=r$. Since $Y_1^L\cup Y_2^L$ is null in $L$ there are uncountably many constructible reals $(x_\alpha:\alpha<\aleph_1)$ not belonging to $Y_1^L\cup Y_2^L$, and therefore not to $Y_1\cup Y_2$.  On the other hand, since $\R^*$ has full measure there is $r\in \R^*\setminus L$. The horizontal section $A_1^{x_\alpha}$ contains only constructible reals since $A_1^{x_\alpha}$ is $\Sigma^1_2(a,x_\alpha)$, and so if it contained a non-constructible real then it would be uncountable by Theorem \ref{t.pst}. Since $A_1\cup A_2$ cover $\R^2$ it must then be the case that the vertical section $(A_2)_r$ contains all the points of the form $(r,x_\alpha)$. But this contradicts that $(A_2)_r$ is countable.
\end{proof}

\subsection{Polarized partitions}

Another type of statement that can be proved by counting arguments analogous to the above are polarized partition relations for $\Sigma^1_2$ colourings of $\mathscr R\times\mathscr R$ (where, as in \S 2, $\mathscr R$ refers to an uncountable recursively presented Polish space.) These may be viewed as regularity properties that $\Sigma^1_2$ colourings have in the presence of a non-constructible real. 

We have the following definable analogue of \cite[24.27]{koto06}:

\begin{theorem}\label{t.partition}
The following are equivalent:
\begin{enumerate}
\item $\mathscr R\nsubseteq L$.
\item For every $\Sigma^1_2$-definable function $f:\mathscr R\times\mathscr R\to\omega$ there are sets $C,D\subseteq\mathscr R$ such that $|C|=|D|=2$ and $f\restrict C\times D$ is monochromatic.
\item For every $\Sigma^1_2$-definable function $f:\mathscr R\times\mathscr R\to\omega$ there are sets $C,D\subseteq\mathscr R$ such that $|C|=|D|=\aleph_0$ and $f\restrict C\times D$ is monochromatic.
\item For every $\Sigma^1_2$-definable function $f:\mathscr R\times\mathscr R\to\omega$ there are countably infinite $\Sigma^1_2$ sets $C,D\subseteq\mathscr R$ such that $f\restrict C\times D$ is monochromatic.
\end{enumerate}
\end{theorem}

\begin{proof}
(4)$\implies$ (3)$\implies $(2) is clear.

(1)$\implies$(3): We may assume that $V=L[z]$ for some $z\notin L$ and that it holds that $\aleph_1^L=\aleph_1$. Assume (3) fails, and fix $f$ witnessing this. For $s\in [\R]^\omega$, let
$$
T(s,i)=\{y\in\R: (\forall x\in s) f(x,y)=i\}.
$$
Since $s$ is a countable sequence this quantification over $s$ may be replaced by a number quantifier over the domain of $s$. Thus $T(s,i)$ is $\Sigma^1_2(s)$. By assumption we have that $|T(s,i)|<\aleph_0$ and so $T(s,i)\subseteq L$ by Theorem \ref{t.pst}. Let $U_i=\{x\in\R\cap L: f(x,z)=i\}$. Then $|U_i|<\aleph_0$ since otherwise we could find $s_i\subseteq U_i$ of size $\aleph_0$ from which $z\in T(s_i,i)$ would follow, contradicting that $z\notin L$. But now we have
$$
\R\cap L=\bigcup_{i\in\omega} U_i
$$ 
so that $|\R\cap L|$ is countable, a contradiction.

(3)$\implies$ (4): Fix a $\Sigma^1_2$-definable function $f$ and $i\in\omega$ such that there exists $C,D\in [\R]^\omega$ with $f(C\times D)=\{i\}$. Since $f$ is $\Sigma^1_2$ the set
$$
\{(C,D)\in [\R]^\omega\times [\R]^\omega: (\forall (x,y)\in C\times D) f(x,y)=i\}
$$
is $\Sigma^1_2$, and it is non-empty by the above. Thus by $\Sigma^1_2$-uniformization (e.g. \cite[4E.4]{moschovakis09}) it contains a $\Sigma^1_2$ definable pair $(C,D)$.

(2)$\implies$(1): Suppose $\R\subseteq L$ and let $\preceq$ denote the usual $\Sigma^1_2$ wellordering of $L\cap\R$. Recall $\IS^\#(x,y)$ from \ref{s.initseg}, and define
$$
f(x,y)=\left\{\begin{array}{ll}
\IS^\#(x,y)+1 & \text{ if } x\prec y\\
\IS^\#(y,x)+1 & \text{ if } y\prec x\\
0 & \text{ if } x=y
\end{array}\right.
$$
Let $\{x,x'\}, \{y,y\}\subseteq\R$ where $x\neq x'$ and $y\neq y'$, and assume that $x,x',y\preceq y'$. Then $f(x,y')\neq f(x',y')$, and so $f\restrict \{x,x'\}\times\{y,y'\}$ is not monochromatic.
\end{proof}

\subsection{A Schur type partition result}

As an application of Theorem \ref{t.partition} we prove the following definable analogue of \cite[24.37]{koto06}.

\begin{theorem}
There is a non-constructible real if and only if for any $\Sigma^1_2$ colouring $g:\R\to\omega$ there are four distinct $x_{00},x_{01},x_{10},x_{11}\in\R$ of the same colour such that
$$
x_{00}+x_{11}=x_{01}+x_{10}.
$$
\end{theorem}
\begin{proof}
Assume $\R\nsubseteq L$ and let $g:\R\to\omega$ be a colouring. By \cite[19.2]{kechris95} we can find a continuous $h:2^\omega\to\R$ such that $h(2^\omega)$ is linearly independent over $\Q$. It may be shown using \cite{kechris72} that this $h$ can be taken to be $\Delta^1_1$. Now let $f:2^\omega\times 2^\omega\to\omega: (x,y)\mapsto g(h(x)+h(y))$. Then by Theorem \ref{t.partition} we can find $x_0\neq x_1$ and $y_0\neq y_1$ such that $f\restrict \{x_0,x_1\}\times\{y_0,y_1\}$ is monochromatic. If we let $x_{ij}=h(x_i)+h(y_j)$ for $0\leq i,j\leq 1$ then clearly $x_{00}+x_{11}=x_{01}+x_{1,0}$ and these are disctinct since $h(2^\omega)$ is linearly independent over $\Q$.

Conversely, assume that $\R\subseteq L$. We define a function $g:\R\to\omega$ by
$$
g(x)=m\iff (\exists\epsilon) R^{\R}(x,\epsilon)\wedge \pi_{\epsilon}(m)=x.
$$
Then $g$ is $\Sigma^1_2$. Let $x_{00},x_{01},x_{10}, x_{11}\in\R$ be distinct, and let $\alpha<\omega$ be least such that $x_{i,j}\in L_\alpha$ for all $0\leq i,j\leq 1$. It cannot be the case that three of the $x_{ij}$ are already in some $L_\beta$ where $\beta<\alpha$, since the $L_\alpha$ are closed under addition. Thus two of the $x_{i,j}$ are in $L_\alpha$ and not in any $L_\beta$ for $\beta<\alpha$. But then these two $x_{i,j}$ are coloured differently by $g$.
\end{proof}

\begin{remark}\label{r.endremark}
It is clear from the above that what is really needed to make all of the above theorems work for $\Sigma^1_n$ (or more generally, $\Sigma^1_n(a)$ versions) is an inner model relative to which we have a $\Sigma^1_n$ absoluteness principle and a perfect set theorem for $\Sigma^1_n$. If we have this, then we will be able to prove that the $\Sigma^1_n$ versions of the statements in Theorem \ref{t.mainthm1} and Theorem \ref{t.partition} are equivalent to all reals being in that inner model.

For example, it is well-known (see \cite[\S 15]{kanamori03}) that if there is a measurable cardinal $\kappa$ and $U$ is an ultrafilter witnessing this, then the inner model $L[U]$ has this relationship to the class of $\Sigma^1_3$ sets, provided that $0^\sharp$ does not exist. Thus in this context we obtain $\Sigma^1_3$ versions of Theorem \ref{t.mainthm1} and Theorem \ref{t.partition}, with $L$ replaced by $L[U]$.

Philip Welch has further pointed out to us that you can more generally do this using
the core model below one Woodin cardinal. Assume (i) there exists a measurable cardinal and (ii) sharps for reals. Then this model is $\Sigma^1_3$ correct, and so the above theorems work over this model.
\end{remark}

\bibliographystyle{amsplain}
\bibliography{CHequivs}

\end{document}